\documentclass[10pt]{amsart}
\newlength{\basicwidth}

\newlength{\shortbasicwidth}

\newlength{\basicheight}\numberwithin{equation}{section}

\begin{document}

\begin{center}
\title
{Inequalities for trigonometric sums}
\maketitle
\end{center}

\vspace{0.4cm}
\begin{center}
HORST ALZER$^a$ \, \, and \,\, MAN KAM KWONG$^b$
\end{center}

\vspace{1.cm}
\begin{center}
$^a$ Morsbacher Stra\ss e 10, 51545 Waldbr\"ol, Germany\\
\emph{Email:} \tt{h.alzer@gmx.de}
\end{center}

\vspace{0.3cm}
\begin{center}
$^b$ Department of Applied Mathematics, The Hong Kong Polytechnic University,\\
Hunghom, Hong Kong\\
\emph{Email:} \tt{mankwong@connect.polyu.hk}
\end{center}

\vspace{1.2cm}
\begin{center}
\emph{To the memory of Richard Bruce Paris}
\end{center}

\vspace{1.2cm}
{\bf{Abstract.}} We present several new  inequalities for trigonometric sums. Among others, we show that the inequality
$$
\sum_{k=1}^n (n-k+1)(n-k+2)k\sin(kx) > \frac{2}{9} \sin(x) \bigl( 1+2\cos(x) \bigr)^2
$$
holds for all $n\geq 1$ and $x\in (0, 2\pi/3)$. The constant factor $2/9$ is sharp. This refines the classical Szeg\"o-Schweitzer inequality which states that the sine sum is positive for all $n\geq 1$ and $x\in (0,2 \pi/3)$. Moreover, as an application of one of our results we obtain a two-parameter class of absolutely monotonic functions.

\vspace{0.8cm}

\vspace{0.2cm}

\newpage

\section{Introduction and statement of the main results}

{\bf{I.}} \, 
In the literature, we can find many papers on inequalities for various trigonometric sums. A reason for the tremendous interest in these inequalities is the fact that they have noteworthy applications, for example, in geometric function theory, number theory, approximation theory and numerical analysis.  Detailed information on this subject with interesting historical comments and many references are given in Askey \cite{A}, Askey and Gasper \cite{AG}, Milovanovi\'c et al. \cite[Chapter 4]{MMV}.
This paper is concerned with some remarkable inequalities for trigonometric sums obtained by the well-known Hungarian mathematicians L. Fej\'er (1880-1959), F. Luk\'acs (1891-1918), G. Szeg\"o (1895-1985), P. Tur\'an (1910-1976) and M. Schweitzer (1923-1945).

\vspace{0.2cm}
{\bf{II.}}  \, In 1935, Tur\'an \cite{T} studied properties of the Ces\`aro means of a sine series. A key role in his investigations plays the elegant inequality
\begin{equation}
\sum_{k=1}^n {n-k+m\choose m}\sin(kx)>0
\end{equation}
which is valid for all natural numbers $m,n$ and real numbers $x\in (0,\pi)$. Extensions, refinements and relatives of (1.1) were given by Alzer and Fuglede \cite{AFu}, Alzer and Kwong \cite{AK1, AK2, AK3}, Bustoz \cite{B}. Our first theorem provides a cosine counterpart of (1.1).

\vspace{0.3cm}
{\bf{Theorem 1.}} \emph{Let $m\geq 1$ be an integer. For all integers $n\geq 1$ and real numbers $x\in (0,\pi)$, we have
\begin{equation}
\sum_{k=0}^n {n-k+m\choose m}\cos(kx)>m.
\end{equation}
The lower bound is sharp.}

\vspace{0.3cm}
{\bf{Remark 1.}} The special case $m=1$ leads to
$$
\sum_{k=0}^n (n-k+1)\cos(kx)>1 \quad(n\geq 1; \, 0<x<\pi).
$$
This is a striking companion to the  Luk\'acs inequality
\begin{equation}
\sum_{k=1}^n (n-k+1)\sin(kx)>0 \quad(n\geq 1; \, 0<x<\pi).
\end{equation}

\vspace{0.3cm}
The following  theorem provides analogues of (1.1) and (1.2).

\vspace{0.5cm}
{\bf{Theorem 2.}} \emph{Let $m\geq 1$ be an integer. For all integers $n\geq 1$ and real numbers $x\in (0,\pi)$, we have
\begin{eqnarray}
\sum_{k=0}^n {n-k+m\choose m}\cos((k+1/2)x)
>
 \left\{ \begin{array}{ll}
-1/4, & \textrm{ {if} $m=1$,}\\
0, & \textrm{ {if} $m\geq 2$,}
\end{array}\right.
\end{eqnarray}
and
\begin{equation}
\sum_{k=0}^n {n-k+m\choose m}\sin((k+1/2)x)>0.
\end{equation}
The given lower bounds are sharp.}

\vspace{0.3cm}
Next, we present inequalities which are closely related to (1.4) and (1.5).

\vspace{0.3cm}
{\bf{Theorem 3.}} \emph{Let $m\geq 1$ be an integer. For all integers $n\geq 1$ and real numbers $x\in (0,\pi)$, we have
\begin{equation}
\sum_{k=0 \atop \, k \,  even}^n {n-k+m\choose m}\cos\bigl(  (k+1/2)x\big)>0
\end{equation}
and
\begin{equation}
\sum_{k=0 \atop \, k \,  even}^n {n-k+m\choose m}\sin\bigl(  (k+1/2)x\big)>0.
\end{equation}
In both cases, the
 lower bound $0$ is sharp.}

\vspace{0.3cm}
{\bf{Remark 2.}} From Theorems 2 and 3 with $m=1$ we obtain the following Luk\'acs-type inequalities which hold for all $n\geq 1$ and $x\in (0,\pi)$,
$$
\sum_{k=0}^n (n-k+1)\cos((k+1/2)x)>-\frac{1}{4}, \quad \, \, 
\sum_{k=0}^n (n-k+1)\sin((k+1/2)x)>0,
$$
$$
\sum_{k=0}^{  [n/2] } (n-2k+1)   \cos\bigl(  (2k+1/2)x\big)>0, 
\quad \, \, 
\sum_{k=0}^{   [n/2]   } (n-2k+1)   \sin\bigl(  (2k+1/2)x\big)>0.
$$
The given lower bounds are sharp.

\vspace{0.3cm}
Tur\'an \cite{T} pointed out that (1.1) can be used to obtain a sine inequality with two variables,
\begin{equation}
\sum_{k=1}^n {n-k+m\choose m}\frac{\sin(kx) \sin(ky)}{k}>0 \quad (m,n\geq 1; \, 0<x,y<\pi).
\end{equation}
We show that an application of (1.5) leads to a counterpart of (1.8).

\vspace{0.3cm}
{\bf{Theorem 4.}}  \emph{Let $m\geq 1$ be an integer. For all integers $n\geq 1$ and real numbers $x,y\in (0,\pi)$, we have
\begin{equation}
\sum_{k=0}^n {n-k+m\choose m}\frac{\sin((2k+1)x) \sin((2k+1)y)}{2k+1}>0.
\end{equation}
The lower bound is sharp.}

\vspace{0.5cm}
{\bf{Remark 3.}}   Inequality (1.9)  with $y=\pi/2$ and $y=\pi/4$ gives
$$
\sum_{k=0}^n (-1)^{\tau_j(k)} {n-k+m\choose m}\frac{\sin((2k+1)x)}{2k+1} >0 \quad (m,n\geq 1; \, 0<x<\pi; \, j=1,2),
$$
where
\begin{displaymath}
\tau_1(k)=k  \quad\mbox{and} \quad \tau_2(k) = \left\{\begin{array}{ll}
k/2,  & \textrm{if $k$ is even,}\\
(k-1)/2, & \textrm{if $k$ is odd.}
\end{array} \right.
\end{displaymath}

\vspace{0.2cm}
{\bf{III.}} \, In 1941, Szeg\"o \cite{S} offered several inequalities for trigonometric sums and used his results to prove the univalence of certain power series. One of his inequalities states that
\begin{equation}
\sum_{k=1}^n (n-k+1)(n-k+2)k\sin(kx)>0\quad (n\geq 1; \, 0<x\leq \theta_0),
\end{equation}
where $\theta_0=1.98231...$. Schweitzer \cite{Sch} refined this result. He proved that the constant $\theta_0$ can be replaced by $2\pi/3$ and that this bound is best possible. Applications and related results can be found in  Alzer and Kwong \cite{AK4, AK5} and Askey and Fitch \cite{AF}.
The next theorem presents a positive minorant for the sine sum in (1.10).

\vspace{0.3cm}
{\bf{Theorem 5.}} \emph{For all $n\geq 1$ and $x\in (0,2\pi/3)$, we have
\begin{equation}
\sum_{k=1}^n (n-k+1)(n-k+2)k\sin(kx) > \lambda \sin(x) \bigl( 1+2\cos(x) \bigr)^2
\end{equation}
with the best possible constant factor $\lambda={2}/{9}$.}

\vspace{0.3cm}
An application of Theorem 5 gives the following cosine inequality.

\vspace{0.3cm}
{\bf{Corollary.}} \emph{For all $n\geq 1$ and $x\in (0,2\pi/3)$, we have
\begin{equation}
\sum_{k=1}^n (n-k+1)(n-k+2)\bigl( 1- \cos(kx) \bigr)
>\mu \bigl( 1- \cos(x) \bigr) \bigl( 13+10\cos(x)+4\cos^2(x)\bigr)
\end{equation}
with the best possible constant factor $\mu={2}/{27}$.}

\vspace{0.3cm}
{\bf{IV.}} \, A function $F:I\rightarrow \mathbb{R}$,  where $I\subset \mathbb{R}$ is an interval, is called absolutely monotonic if $F$ has derivatives of all orders and satisfies
$$
F^{(n)}(x) \geq 0 \quad (n=0,1,2,...; \; x\in I).
$$
These functions play a role in various fields, like for example, the theory of analytic functions and probability theory. We refer to Boas \cite{Boas} and Widder \cite[Chapter IV]{Wi} for more information on this subject. It is known that inequalities for trigonometric sums can be applied to prove that certain functions are absolutely monotonic; see Milovanovi\'c et al. \cite[Chapter 4.2.5]{MMV}.
Here, we use Theorem 1 to present a new two-parameter class of absolutely monotonic functions.

\vspace{0.3cm}
{\bf{Theorem 6.}} \emph{Let $m\geq 1$ be an integer and let $\omega\in [-1,1]$. The function
$$
W_{m,\omega}(x)=m-1-\frac{m}{1-x}+\frac{1-\omega x}{  (1-x)^{m+1} (1-2 \omega x+x^2)}
$$
is absolutely monotonic on $(0,1)$.}

\vspace{0.3cm}
{\bf{Remark 4.}}  Applying  Theorem 6 and  the Petrovi\'c functional inequality for convex functions (see Mitrinovi\'c \cite[Section 1.4.7]{M}) gives that $W_{m,\omega}$ is superadditive on $(0,1)$. This means that if $m\geq 1$ and $\omega \in [-1,1]$, then we have for all nonnegative real numbers $x,y$ with $x+y<1$,
$$    
W_{m,\omega}(x)+W_{m,\omega}(y)\leq W_{m,\omega}(x+y).
$$

\vspace{0.2cm}
{\bf{V.}}  \, In the next section, we collect some helpful lemmas. The proofs of the theorems and the corollary are given in Section 3 to 
Section 8. The numerical and algebraic computations have been carried out using the computer software Maple 13.

\section{Lemmas}

The first two lemmas present inequalities for certain classes of cosine and sine sums. Both results are due to Fej\'er \cite{F1, F2}.

\vspace{0.3cm}
{\bf{Lemma 1.}} \emph{Let $c_0, c_1, ..., c_N$ be real numbers such that
\begin{equation}
c_0-c_1 \geq c_1 -c_2 \geq \cdots \geq c_{N-1}-c_N \geq c_N\geq 0.
\end{equation}
Then, for $x\in \mathbb{R}$,}
$$
\frac{c_0}{2}+\sum_{k=1}^N c_k \cos(kx)\geq 0.
$$

\vspace{0.3cm}
{\bf{Lemma 2.}} \emph{Let $c_1,...,c_N$ be nonnegative real numbers. If
$$
\sum_{k=1}^N k c_k \sin(kt)>0\quad(0<t<\pi),
$$
then}
$$
\sum_{k=1}^N c_k \sin(kx)\sin(ky)>0 \quad (0<x,y<\pi).
$$

\vspace{0.3cm}
The following lemmas are needed in the proof of Theorem 5. First, we collect some properties of the functions
\begin{equation}
S_n(x)=(18n+24)\sin(x) -(9n+27)\sin((n+1)x) +9n\sin((n+2)x) +2\sin(4x)-\sin(5x)
\end{equation}
and
$$
 L_n(x)=(18n+24)\sin(x)-18n \sin(x/2)-29.1.
$$

\vspace{0.3cm}
{\bf{Lemma 3.}} \emph{Let $n\geq 21$ and $x\in (0,2\pi/3)$. Then $S_n(x)> L_n(x)$.}

\vspace{0.3cm}
\begin{proof}
We have
$$
2\sin(4x)-\sin(5x)+2.1\sin(x)=\sin(x) v(\cos(x))
$$
with
$$
v(t)=-16t^4+16t^3+12t^2-8t+1.1.
$$
Next, we apply Sturm's theorem to determine the number of distinct roots of an algebraic polynomial located in an interval; see van der Waerden \cite[Section 79]{W}.
We obtain that $v$ has no zero on $[-1/2,1]$, so that $v(1)=5.1$ leads to $v(t)>0$ for $t\in [-1/2,1]$. This gives
$$
2\sin(4x)-\sin(5x)>-2.1\sin(x)\geq -2.1.
$$
It follows that
\begin{eqnarray}\nonumber
S_n(x) & = & (18n+24)\sin(x)+18n\sin(x/2)\cos((n+3/2)x) -27\sin((n+1)x)+2\sin(4x)-\sin(5x) \\ \nonumber
& > & (18n+24)\sin(x)-18n\sin(x/2)-27-2.1=L_n(x). \nonumber
\end{eqnarray}.
\end{proof}

\vspace{0.3cm}
{\bf{Lemma 4.}} \emph{Let $n\geq 21$. The function $L_n$ is concave on $(0,2\pi/3)$.}

\vspace{0.3cm}
\begin{proof}
We have
$$
-L''_n(x)=18n\sin(x/2)\Bigl( 2\cos(x/2)-\frac{1}{4}\Bigr)+24\sin(x)\geq \frac{27}{2} n \sin(x/2)+24\sin(x)>0.
$$
\end{proof}

\vspace{0.3cm}
{\bf{Lemma 5.}} \emph{We have}
$$
L_n\Bigl( \frac{1.1 \pi}{n}\Bigr)>0 \quad (n\geq 2) \quad{and} \quad L_n\Bigl( \frac{2\pi}{3}-\frac{1}{n}\Bigr)>0 \quad (n\geq 21).
$$

\vspace{0.3cm}
\begin{proof}
(i) We have $L_2(1.1\pi/2)=2.78...$. Let $n\geq 3$. Using
\begin{equation}
x-\frac{1}{6}x^3\leq \sin(x)\leq x \quad (x\geq 0)
\end{equation}
gives
\begin{eqnarray}\nonumber
L_n\Bigl( \frac{1.1\pi}{n}\Bigr) & = & (18n+24)\sin \Bigl( \frac{1.1\pi}{n}\Bigr)-18n \sin\Bigl( \frac{1.1\pi}{2n}\Bigl) -29.1\\ \nonumber
& \geq & (18n+24) \Bigl(  \frac{1.1\pi}{n}-\frac{1}{6}\Bigl(  \frac{1.1\pi}{n}\Bigr)^3 \Bigr) -18n \cdot \frac{1.1\pi}{2n}-29.1 \\ \nonumber
& = & \frac{Y(n)}{n^3}\nonumber
\end{eqnarray}
with
$$
Y(n)=an^3+bn^2-cn-d,
$$
$$
a=9.9\pi-29.1=2.00..., \quad b=26.4\pi=82.93..., \quad c=3.993 \pi^3 =123.80..., \quad d=5.324 \pi^3 =165.07... .
$$
Since $Y$ is positive on $[3,\infty)$, we conclude that $L_n(1.1\pi/n)>0$.

(ii) Let $n\geq 21$. We have
\begin{eqnarray}\nonumber
L_n\Bigl( \frac{2\pi}{3}-\frac{1}{n}\Bigr) & = & (18n+24) \sin\Bigl(\frac{\pi}{3}+\frac{1}{n}\Bigr)-18n\cos\Bigl( \frac{\pi}{6}+\frac {1}{2n}\Bigr)-29.1\\ \nonumber
& \geq &  (18n+24) \sin\Bigl(\frac{\pi}{3}+\frac{1}{n}\Bigr) -18n\cos(\pi/6)-29.1 \\ \nonumber
& = & 9n\sin(1/n)+9\sqrt{3} n \cos(1/n) -9\sqrt{3}n +24\sin\Bigl(\frac{\pi}{3}+\frac{1}{n}\Bigr)-29.1 \\ \nonumber
& \geq & 9n\sin(1/n)-9\sqrt{3}n\bigl( 1-\cos(1/n)\bigr)+12\sqrt{3}-29.1.\nonumber
\end{eqnarray}
We set
$$
\alpha=189\sin(1/21)+12\sqrt{3}-29.1.
$$
Since $x\mapsto \sin(x)/x$ is decreasing on $(0,\pi]$ and
$$
1-\frac{x^2}{2}\leq \cos(x) \quad (x\geq 0),
$$
we obtain
$$
L_n\Bigl( \frac{2\pi}{3}-\frac{1}{n}\Bigr)\geq \alpha -9\sqrt{3}n \bigl(1-\cos(1/n)\bigr)\geq \alpha-\frac{9\sqrt{3}}{42}=0.31... .
$$
\end{proof}

\vspace{0.3cm}
{\bf{Lemma 6.}} \emph{Let $n=3m$ with $m\geq 7$ and $x\in (2\pi/3-1/n, 2 \pi/3)$. Then}
$$
S''_n(x)=-(18n+24)\sin(x)+(9n+27)(n+1)^2\sin((n+1)x) -9n(n+2)^2 \sin((n+2)x) -32\sin(4x)+25\sin(5x)>0.
$$

\vspace{0.3cm}
\begin{proof}
We have
$$
\sin((n+1)x) \geq \sin\Bigl(  \frac{2\pi}{3}-\frac{22}{21}\Bigr)=0.865...
\quad\mbox{and}
\quad \sin((n+2)x)\leq \sin\Bigl(  \frac{4\pi}{3}-\frac{23}{21}\Bigr)=0.048... .
$$
Thus
\begin{eqnarray}\nonumber
S''_n(x) &   \geq & -(18n+24)\sin(x)+0.86 (9n+27)(n+1)^2 -9\cdot 0.05 n (n+2)^2 -32\sin(4x)+25\sin(5x)\\ \nonumber
& \geq & -(18n+24)+0.86 (9n+27) (n+1)^2 -0.45 n(n+2)^2-57\\ \nonumber
& = & 7.29n^3+36.9 n^2 +34.38 n -57.78>0. \nonumber
\end{eqnarray}
\end{proof}

\vspace{0.3cm}
Moreover, we need lower bounds for the functions
\begin{equation}
f_n(t)=24\sin\Bigl(\frac{t}{n+2}\Bigl)+2\sin\Bigl(\frac{4t}{n+2}\Bigl)-\sin\Bigl(\frac{5t}{n+2}\Bigl),
\end{equation}
\begin{equation}
g_n(t)=18n\sin\Bigl(\frac{t}{n+2}\Bigl)-27\sin\Bigl(\frac{(n+1)t}{n+2}\Bigl),
\end{equation}
\begin{equation}
h_n(t)=9n\sin(t)-9n\sin\Bigl(\frac{(n+1)t}{n+2}\Bigl)
=18n\sin\Bigl(\frac{t}{2n+4}\Bigl)
\cos\Bigl(\frac{(2n+3)t}{2n+4}\Bigl).
\end{equation}

{\bf{Lemma 7.}} \emph{Let $n\geq 21$ and $t\in (2.5,1.21\pi)$. Then}
$$
f_n(t)\geq \frac{26.3 t}{n+2}, \quad g_n(t)>9t, \quad h_n(t)\geq -9t.
$$

\vspace{0.3cm}
\begin{proof}
(i) Using (2.3) gives
\begin{eqnarray}\nonumber
f_n(t) & \geq & 24 \Bigl(  \frac{t}{n+2}-\frac{1}{6}\Bigl( \frac{t}{n+2}\Bigr)^3\Bigr)+2\Bigl( \frac{4t}{n+2}-\frac{1}{6}\Bigl(  \frac{4t}{n+2}\Bigr)^3\Bigr)-\frac{5t}{n+2} \\ \nonumber
& = & \frac{t}{n+2}\Bigl( 27-\frac{76}{3}\Bigl( \frac{t}{n+2}\Bigr)^2 \Bigl) \\ \nonumber
& \geq & \frac{t}{n+2}\Bigl(  27-\frac{76}{3}\Bigl(  \frac{1.21\pi}{23}\Bigr)^2\Bigr) \\ \nonumber
& \geq & \frac{26.3 t}{n+2}.\nonumber
\end{eqnarray}

(ii)  Since the sequences $n\mapsto n\sin(t/(n+2))$ and $n\mapsto -\sin((n+1)t/(n+2))$ are increasing, we conclude that $n\mapsto g_n(t)$ is increasing. It follows that
$g_n(t)\geq g_{21}(t)$. 
Let
$$
G(t)=\frac{1}{27}\Bigl( g_{21}(t)-9t\Bigr)=14\sin\Bigl(  \frac{t}{23}\Bigr)
-\sin\Bigl(  \frac{22t}{23}\Bigr)-\frac{t}{3}
$$
The functions
$t\mapsto \sin (t/23)$ and $t\mapsto -\sin(22t/23)$ are increasing on $[2.5,1.21\pi]$. Let $2.5\leq r\leq t\leq s\leq 1.21 \pi$. Then
we obtain
$$
G(t)\geq 14\sin\Bigl(  \frac{r}{23}\Bigr)
-\sin\Bigl(  \frac{22r}{23}\Bigr)-\frac{s}{3}=H(r,s), \quad\mbox{say}.
$$
By direct computation, we get
$$
H\Bigl( 2.5+\frac{k}{100}, 2.5+\frac{k+1}{100}\Bigr)>0 \quad (k=0,1,...,39),
\quad H(2.9, 1.21\pi)=0.13... .
$$
It follows that $G(t)>0$ for $t\in (2.5,1.21\pi)$. This leads to $g_n(t)>9t$.

(iii) We have
$$
h_n(t)\geq -18n \sin\Bigl(  \frac{t}{2n+4}\Bigr)\geq -18n \cdot  \frac{t}{2n+4}\geq -9t.
$$
\end{proof}

\vspace{0.3cm}
\section{Proof of Theorem 1}

Let $m,n\geq 1$. We set
$$
c_k={n-k+m\choose m}\quad (k=0,1,...,n+1).
$$
Then,
$$
c_{k+2}-2c_{k+1}+c_k=\frac{m(m-1)}{(n-k)(n-k+m-1)}c_{k+1}\geq 0 \quad (k=0,...,n-1).
$$
It follows that (2.1) holds, so that we obtain, for $x\in\mathbb{R}$,
\begin{equation}
T_n(m,x)=\frac{1}{2}{n+m\choose m}+\sum_{k=1}^n {n-k+m\choose m}\cos(kx)\geq 0.
\end{equation}
We denote the cosine sum in (1.2) by $B_n(m,x)$. Let $x\in (0,\pi)$. Then
$$
B_1(m,x)-m=1+\cos(x)>0.
$$
Let $n\geq 2$. We obtain
\begin{equation}
B_n(m,x)-m-T_n(m,x)=\frac{1}{2}{n+m\choose m}-m\geq \frac{1}{4}\bigl(  (m-1)m+2\bigr)>0.
\end{equation}
From (3.1) and (3.2) we conclude that $B_n(m,x)>m$. Since $B_1(m,\pi)=m$, it follows that $m$ is the best possible lower bound in (1.2).

\vspace{0.5cm}
\section{Proof of Theorem 2}

\emph{Proof of} (1.4).   We denote the sum in (1.4) by $U_n(m,x)$. Then we have
$$
U_n(1,x)=\frac{\cos(x/2)-\cos((n+3/2)x)}{2(1-\cos(x))}\geq \frac{\cos(x/2)-1}{2(1-\cos(x))}=-\frac{1}{4(1+\cos(x/2))}>-\frac{1}{4}.
$$
This settles (1.4) for $m=1$. Moreover, if we set $x_n=4n\pi/(4n+1)$, then
$$
U_{2n-1}(1,x_n)=\frac{\cos(x_n/2)-1}{2(1-\cos(x_n))}.
$$
Since
$$
\lim_{n\to\infty} U_{2n-1}(1,x_n)=-\frac{1}{4},
$$
we conclude that the lower bound $-1/4$ is sharp.

Next, let $m\geq 2$.  Then
\begin{equation}
U_1(m,x)=(m+1)\cos(x/2)+\cos(3x/2)\geq 3\cos(x/2)+\cos(3x/2)=4\cos^3(x/2)>0.
\end{equation}
We have
$$
U_n(2,x)=\frac{1}{2}\sum_{k=0}^n (n-k+1)(n-k+2)\cos((k+1/2)x)=\frac{1}{2\sin(x/2)}\sum_{k=0}^n (n-k+1)\sin((k+1)x),
$$
so that (1.3) yields
\begin{equation}
U_n(2,x)>0.
\end{equation}
Since
$$
{N+1\choose \nu}={N\choose \nu}+{N\choose \nu-1},
$$
we obtain the representation
\begin{equation}
U_{n+1}(m+1,x)=U_{n+1}(m,x)+U_n(m+1,x).
\end{equation}
Using (4.1), (4.2) and (4.3) we obtain by induction that $U_n(m,x)>0$ for all $n\geq 1$ and $m\geq 2$. 
Since $U_n(m,\pi)=0$, we conclude that the lower bound $0$ is best possible.

\vspace{0.3cm}
\emph{Proof of} (1.5). We denote the  sum in  (1.5) by $V_n(m,x)$. Then
$$
V_n(m,x)=\cos(x/2) A_n(m,x) + \sin(x/2) B_n(m,x),
$$
where $A_n(m,x)$ and $B_n(m,x)$ are the sums given in (1.1) and (1.2), respectively. 
Using (1.1) and (1.2) gives 
$$
V_n(m,x)>0+\sin(x/2)\cdot m >0.
$$
Moreover, since $V_n(m,0)=0$, it follows that the lower bound $0$ is best possible.

\vspace{0.5cm}
\section{Proof of Theorem 3}

We denote  the sums in (1.6) and (1.7) by $C_n(m,x)$ and $D_n(m,x)$, respectively. 
Since $C_n(m, \pi-x)=D_n(m,x)$, it suffices to prove that $D_n(m,x)>0$. First, we consider the case $m=1$. 
We have
$$
\frac{1}{2} D_{2n+1}(1,x)=\sum_{k=0}^n (n-k+1)\sin( (2k+1/2)x) =\frac{E_n(x)}{32\sin^3(x/2) \cos^2(x/2)}
$$
with
$$
E_n(x)=\sin(x) \bigl(  2(n+1)\sin(x)-\sin( 2(n+1)x ) \bigr) +4\sin^2(x/2) \sin^2((n+1)x).
$$
Since
$$
N\sin(x) >\sin(Nx)  \quad (N=2,3,....; \, 0<x<\pi),
$$
we conclude that $E_n(x)>0$. It follows that 
\begin{equation}
D_{2n+1}(1,x)>0 \quad (n\geq 0).
\end{equation}
Let $r\in (-1,1)$. We define
$$
J_x(r)=\sum_{k=0}^\infty \frac{\sin((2k+1/2)x)}{\sin(x)}r^k,
\quad 
K_x(r)=\sin(x)\frac{1+r}{2(1-r)^2} J_x(r),
\quad
M_x(r)=\sin(x)\frac{1}{(1-r)^2} J_x(r).
$$
Using
$$
\sum_{n=0}^\infty (n+1)r^n=\frac{1}{(1-r)^2}\quad\mbox{and}\quad \sum_{n=0}^\infty (n+1/2)r^n =\frac{1+r}{2(1-r)^2}
$$
gives
$$
K_x(r)=\sum_{n=0}^\infty \sum_{k=0}^n (n-k+1/2)\sin((2k+1/2)x) r^n=\frac{1}{2}\sum_{n=0}^\infty D_{2n}(1,x) r^n
$$
and
$$
M_x(r)=\sum_{n=0}^\infty \sum_{k=0}^n (n-k+1)\sin((2k+1/2)x) r^n=\frac{1}{2}\sum_{n=0}^\infty D_{2n+1}(1,x) r^n.
$$
Since
$$
2K_x(r)=(1+r) M_x(r),
$$
we obtain
\begin{equation}
 D_{2n}(1,x)=\frac{1}{2}( D_{2n-1}(1,x)+D_{2n+1}(1,x)) \quad (n\geq 1).
\end{equation}
From (5.1) and (5.2) we conclude that $D_{2n}(1,x)>0$ $(n\geq 1)$.

Next, let
 $m\geq 2$ and let $U_n(m,x)$ and $V_n(m,x)$ be the sums given in (1.4) and (1.5), respectively.
Applying Theorem 2 gives
$$
D_n(m,x)
=\frac{1}{2} \bigl( U_n(m,\pi-x)+V_n(m,x) \bigr) >0.
$$
Since $C_n(m,\pi)=D_n(m,0)=0$, we conclude  that $0$ is the best possible lower bound in (1.6) and (1.7).

\vspace{0.5cm}
\section{Proof of Theorem 4}

Let $V_n(m,x)$ be the sine sum in (1.5). 
Since
$$
V_n(m, \pi+x)=V_n(m, \pi-x)  \quad\mbox{and}\quad V_n(m,\pi)>0,
$$
we conclude from  Theorem 2 that
$$
V_n(m,2t)>0 \quad (0<t<\pi).
$$
We set
$$
c_{2k}=0, \quad
c_{2k+1}=\frac{1}{2k+1}{n-k+m\choose m} \quad (k=0,1,...,n).
$$
Then we have
$$
V_n(m,2t)=\sum_{k=1}^{2n+1} k c_k \sin(kt)>0,
$$
so that  Lemma 2 with $N=2n+1$ gives for $x,y\in (0,\pi)$,
$$
\sum_{k=1}^{2n+1}c_k \sin(kx)\sin(ky)=\sum_{k=0}^n {n-k+m\choose m}\frac{\sin((2k+1)x)\sin((2k+1)y)}{2k+1}>0.
$$
If we set $x=0$, then equality holds in (1.9). This implies that the lower bound $0$ is sharp.

\vspace{0.5cm}
\section{Proof of Theorem 5 and the Corollary}

\emph{Proof of Theorem} 5.
Let $P_n(x)$ be the sum in (1.11). We define
\begin{equation}
Q_n(x)=\frac{1}{\sin(x)}\Bigl(  P_n(x)-\frac{2}{9} \sin(x)  ( 1+2\cos(x) )^2\Bigr).
\end{equation}
Then, with $t=\cos(x)  \in (-1/2,1)$,
$$
Q_1(x)=\frac{8}{9}(2+t)(1-t)>0,
$$
$$
Q_2(x)=\frac{4}{9} (13+16t-2t^2)>0,
$$
$$
Q_3(x)=\frac{52}{9} (1+2t)^2>0.
$$
Let $n\in \{4,5,...,20\}$. Then we have
\begin{equation}
Q_n(x)=R_n(t), \quad t=\cos(x),
\end{equation}
where $R_n$ is an algebraic polynomial of degree $n-1$. Applying Sturm's theorem gives that if $n \not\equiv 0 \, (\mbox{mod} \, 3)$, then $R_n$ has no zero on $[-1/2,1]$, and 
if $n \equiv 0 \, (\mbox{mod} \, 3)$, then $R_n$ has precisely one zero on $[-1/2,1]$, namely, $t=-1/2$. Since $R_n(1)>0$, we conclude that $R_n$ is positive on $(-1/2,1)$. From (7.1) and (7.2) we conclude that  (1.11) holds.

Let $n\geq 21$. First, we prove that (1.11) is valid for $x\in (0, 2.5/(n+2)]$. Using
$$
\frac{2}{9} \sin(x)  \bigl(1+2\cos(x)\bigr)^2=\frac{4}{9} \sin(x)+\frac{4}{9} \sin(2x)+\frac{2}{9}\sin(3x)
$$
gives
$$
P_n(x)-\frac{2}{9}\sin(x) \bigl( 1+2\cos(x)\bigr)^2 =\sum_{k=1}^n a_{k,n} \sin(kx)
$$
with
\begin{displaymath}
a_{k,n} = \left\{ \begin{array}{ll}
n(n+1)-4/9, \, \,  \textrm{if $k=1$,}\\
2(n-1)n -4/9, \, \, \textrm{if $k=2$,}\\
3(n-2)(n-1)-2/9, \, \,  \textrm{if $k=3$,}\\
(n-k+1)(n-k+2) k, \, \,  \textrm{if $k\geq 4$.}
\end{array} \right.
\end{displaymath}
Since $a_{k,n}>0$  and $0<kx<\pi$ $(1\leq k\leq n)$, we conclude that (1.11) holds.

We have
$$
72 \sin^4(x/2) \Bigl( P_n(x)-\frac{2}{9} \sin(x)\bigl( 1+2\cos(x)  \bigr)^2  \Bigr)=S_n(x),
$$
where $S_n(x)$ is given in (2.2).
It remains to prove 
that $S_n$ is positive on $(2.5/(n+2),  2\pi/3)$.  We consider three cases.

\vspace{0.2cm}
{\underline{Case 1.}} $x\in (2.5/(n+2), 1.1\pi/n)$.\\
We set  $t=(n+2)x$. Then  $t \in (2.5,1.21 \pi)$.
Since
$$
S_n\Bigl( \frac{t}{n+2}\Bigr)=f_n(t) +g_n(t) + h_n(t),
$$
where $f_n$, $g_n$ and $h_n$ are defined in (2.4), (2.5) and (2.6), we conclude from Lemma 7  that $S_n(x)>0$.

\vspace{0.2cm}
{\underline{Case 2.}} $x\in [1.1\pi/n, 2\pi/3-1/n]$.\\
Applying Lemmas 3, 4 and 5 yields
$$
S_n(x) > L_n(x) \geq \min \Bigl\{ L_n\Bigl( \frac{1.1\pi}{n}\Bigr),  L_n \Bigl(  \frac{2\pi}{3}-\frac{1}{n}\Bigr) \Bigr\}>0.
$$

\vspace{0.2cm}
{\underline{Case 3.}} $ x\in (2\pi/3-1/n,  2\pi/3)$.\\
We consider three subcases.\\
{\underline{Case 3.1.}} $n=3m+1$.\\
Using
$$
\sin(x)\geq \sin(2\pi/3)=\frac{\sqrt{3}}{2}, \quad \sin((n+1)x)\leq \sin(4\pi/3-22/21)=0.0004...
$$
leads to
$$
S_n(x) \geq  9(\sqrt{3}-1.0005)n +12\sqrt{3}-3.0135>0.
$$

\vspace{0.2cm}
{\underline{Case 3.2.}} $n=3m+2$.\\
Since
$$
\sin(x)\geq \frac{\sqrt{3}}{2},  \quad \sin((n+1)x)\leq 0, \quad \sin((n+2)x) \geq 0,
$$
we obtain
$$
S_n(x)\geq  (9n+12)\sqrt{3}-3>0.
$$

\vspace{0.2cm}
{\underline{Case 3.3.}} $n=3m$.\\
We obtain
$$
S_n(2\pi/3)=S'_n(2\pi/3)=0,
$$
so that Lemma 6 gives $S_n(x)>0$.
This completes the proof of (1.11). 

\vspace{0.2cm}
Since
$$
\lim_{x\to 0}  \frac{P_1(x)}{\sin(x)  (1+2\cos(x))^2}=\lim_{x\to 0}\frac{2}{(1+2\cos(x))^2}=
\frac{2}{9},
$$
we conclude that the constant factor $2/9$ is best possible.

\vspace{0.5cm}
\emph{Proof of the Corollary.}
We denote the cosine sum in (1.12) by $\Theta_n(x)$. From Theorem 5 we obtain
$$
\Theta_n(x)=\int_0^x P_n(t)dt >\frac{2}{9} \int_0^x \sin(t) \bigl( 1+2\cos(t) \bigr)^2 dt=\frac{2}{27}
 \bigl( 1- \cos(x) \bigr) \bigl( 13+10\cos(x)+4\cos^2(x)\bigr).
$$
This settles (1.12). Moreover, since
$$
\lim_{x\to 0}\frac{\Theta_1(x)}{ ( 1- \cos(x) ) ( 13+10\cos(x)+4\cos^2(x)  )    }=\frac{2}{27},
$$
we conclude that $2/27$ is the best possible constant factor in (1.12).

\vspace{0.5cm}
\section{Proof of Theorem 6}

Let $m\geq 1$,  $\theta\in [0,\pi]$ and $x\in (0,1)$. We have
$$
\Lambda_m(x)=\sum_{n=0}^\infty {n+m\choose m} x^n=\frac{1}{(1-x)^{m+1}}
\quad \mbox{and}\quad
\Phi_{\theta}(x)=\sum_{n=0}^\infty \cos(n\theta) x^n=\frac{1-\cos(\theta) x}{1-2 \cos(\theta) x+x^2}.
$$
The Cauchy product formula yields
\begin{eqnarray}\nonumber
\Lambda_m(x) \Phi_{\theta}(x)-\frac{m}{1-x} & = & \sum_{n=0}^\infty \Bigl( \sum_{k=0}^n   {n-k+m\choose m}     \cos(k\theta)-m  \Bigr) x^n \\ \nonumber
& = & 1-m+ \sum_{n=1}^\infty \Bigl( \sum_{k=0}^n  {n-k+m\choose m}           \cos(k\theta)-m  \Bigr) x^n. \\ \nonumber
\end{eqnarray}
Using (1.2) gives that the function
$$
x\mapsto m-1-\frac{m}{1-x} +\Lambda_m(x) \Phi_{\theta}(x)=m-1-\frac{m}{1-x}+\frac{1-\cos(\theta) x}{(1-x)^{m+1} (1-2\cos(\theta) x+x^2)}
$$
is absolutely monotonic on $(0,1)$.

\end{document}